\documentclass[11pt]{article}
\usepackage{amsthm,amsmath,latexsym,amssymb}
\newcommand{\bP}{{\rm |\kern-.15em P}}
\newcommand{\Q}{\kern.3em\rule{.07em}{.65em}\kern-.3em{\rm Q}}
\newcommand{\R}{{\rm I\kern-.15em R}}
\newcommand{\D}{{\rm |\kern-.15em D}}
\newcommand{\h}{{\rm |\kern-.15em H}}
\newcommand{\C}{\kern.3em\rule{.07em}{.65em}\kern-.3em{\rm C}}
\newcommand{\T}{{\rm T\kern-.35em T}}

\theoremstyle{plain}
\newtheorem{theorem}{Theorem}[section]
\newtheorem{lemma}[theorem]{Lemma}

\newtheorem{proposition}[theorem]{Proposition}

\theoremstyle{definition}
\newtheorem{definition}[theorem]{Definition}
\newtheorem{example}[theorem]{Example}

\theoremstyle{remark}
\newtheorem{remark}[theorem]{Remark}

\newcommand\blfootnote[1]{%
  \begingroup
  \renewcommand\thefootnote{}\footnote{#1}%
  \addtocounter{footnote}{-1}%
  \endgroup
}

\begin{document}
\title{Extreme points and support points of conformal mappings}
\author{Ronen Peretz}
 
\maketitle

\begin{abstract}
There are three types of results in this paper. The first, extending a representation theorem on
a conformal mapping that omits two values of equal modulus. This was due to Brickman and Wilken.
They constructed a representation as a convex combination with two terms. Our representation constructs
convex combinations with unlimited number of terms. In the limit one can think of it as an integration
over a probability space with the uniform distribution.
The second result determines the sign of $\Re L(\overline{z}_0(f(z))^2)$ up to a remainder term which is
expressed using a certain integral that involves the L\"owner chain induced by $f(z)$, for a support point $f(z)$
which maximizes $\Re L$. Here $L$ is a continuous linear functional on $H(U)$, the topological vector space of the holomorphic functions
in the unit disk $U=\{z\in\mathbb{C}\,|\,|z|<1\}$. Such a support point is known to be a slit mapping and $f(z_0)$
is the tip of the slit $\mathbb{C}-f(U)$.
The third demonstrates some properties of support points of the subspace $S_n$ of $S$. $S_n$ contains
all the polynomials in $S$ of degree $n$ or less. For instance such a support point $p(z)$ has a zero
of its derivative $p'(z)$ on $\partial U$.
\end{abstract}

\section{Introduction}\label{sec1}

\blfootnote{\textup{2010} \textit{Mathematics Subject Classification}: \textup{30C20, 30C50, 30C55, 30C70, 30C75,
46A03, 46A55}}
\blfootnote{\textit{Key Words and Phrases:} \textup{extreme points, support points, conformal mappings, 
schlicht functions, }}

Let $S:=\{f\in H(U)\,|\,f(0)=f'(0)-1=0, f\,\,{\rm is}\,\,{\rm injective}\,\,{\rm on}\,\,U:=\{z\in\mathbb{C}\,|\,
|z|<1\}\}$. This is the family of normalized conformal mappings on the open unit disk $U$. $S$ is a normal
family and a compact subspace of the holomorphic functions on $U$, $H(U)$. The topology is taken to be that
of uniform convergence on compact subsets of $U$. This topology is locally convex on $H(U)$. We recall
the following standard definitions.

\begin{definition}\label{def1}
Let $X$ be a topological vector space over the field of complex numbers. Let $Y$ be a subset of $X$. 
A point $x\in Y$ is called an extreme point of $Y$ if it has no representation of the form
$$
x=t\cdot y+(1-t)\cdot z,\,\,\,\,\,0<t<1,
$$
as a proper convex combination of two distinct points $y$ and $z$ in $Y$. A point $x\in Y$ is called a support
point of $Y$ if there is a continuous linear functional $L$ on $X$, not constant on $Y$, such that
$$
\Re\{L(x)\}\ge\Re\{L(y)\}\,\,\,{\rm for}\,\,{\rm all}\,\,y\in Y.
$$
\end{definition}
In this paper we will give an extension of a result of L. Brickman and D. R. Wilken. This result whose
elegant proof is essentially due to Brickman and Wilken can be found in \cite{b}. See also \cite{bw}.

Another property we will prove is that for a function $f\in S$ that maximizes $\Re\{L(g)\}$, $g\in S$ where
$L$ is a linear continuous functional on $H(U)$, we have for any natural number $n\in\mathbb{Z}^+$, and for
any positive real number $t\in\mathbb{R}^+$:
$$
\Re L\left\{\overline{z}_0f(z)^2e^{-t}\right\}+
\Re L\left\{\int_t^{\infty}\left\{\frac{e^s f(z,s)(k(s)f(z,s))^2}{1-k(s)f(z,s)}\right\}ds\right\}+o(e^{-t})\le 0.
$$
where $f(z_0)$ is the tip of the monotone slit $\mathbb{C}-f(U)$, $|z_0|=1$ and $f'(z_0)=0$. $f(z,s)$, ($z\in U$,
$s\in\mathbb{R^+}$) is the L\"owner chain generated by the support point $f(z)$. We will use
as a general reference the book \cite{d}. Especially Chapter 9, 275-287 and Chapter 3, 76-113. \\
In the final section we will prove that properties of the support points $f$ of $S$, such as that
$f'$ has a zero on the boundary $\partial U$, are inherited by much smaller subfamilies of $S$ such as $S_n$, the
spaces of all the polynomials in $S$ of degree $n$ or less ($n\in\mathbb{Z}^+$). Clearly the $S_n$'s are less geometric
than $S$. Nevertheless the birth of the slit structure of the image is starting to be visible by their support points.
An important part of geometric function theory is the solution of extremal problems, such as
coefficient problems, integral means problems, distortion problems and many other extremal problems. In order to apply
functional analytic tools it is natural to identify the extreme points of $S$ and its support points. By the Krein-Milman
theorem, there is an extreme point of $S$ among the support points associated with each linear continuous functional on
$H(U)$. Knowing properties of support points might allow restricting the search for a solution to a much smaller
family of points in $S$, than the whole of $S$. This is one aspect of the importance of such results.

\section{A simple extension of a result of Brickman and Wilken}\label{sec2}

Here is a result of Brickman and Wilken, \cite{b}. \\
\\
{\bf Theorem (Brickman and Wilken, \cite{b}).} {\it If a function $f\in S$ omits two values of equal modulus,
then $f$ has the form $f=t\cdot f_1+(1-t)\cdot f_2$, $0<t<1$, where $f_1$ and $f_2$ are distinct functions in $S$
which omit open sets.} \\
\\
The clever proof given by Brickman and Wilken considers the image of $f$, $D=f(U)$ which omits $\alpha$ and $\beta$,
$\alpha\ne\beta$. They define an analytic single-valued branch of $\Psi(w)=\{(w-\alpha)(w-\beta)\}^{1/2}$ in $D$ and
prove that the two functions $w\pm\Psi(w)$ are univalent and have disjoint images of $D$. They normalize to
get two conformal mappings later on that belong to $S$
$$
\Psi_1(w)=\frac{w+\Psi(w)-\Psi(0)}{1+\Psi'(0)},\,\,\,\Psi_2(w)=\frac{w-\Psi(w)+\Psi(0)}{1-\Psi'(0)}.
$$
Now, by the identity
$$
(1+\Psi'(0))\cdot\Psi_1(w)+(1-\Psi'(0))\cdot\Psi_2(w)=2\cdot w,
$$
and with the compositions $f_1=\Psi_1\circ f$, $f_2=\Psi_2\circ f$ they obtain two functions $f_1$ and $f_2$
in $S$ that satisfy $f(z)=t\cdot f_1(z)+(1-t)\cdot f_2(z)$ for $z\in U$, where 
$$
t=\frac{1}{2}(1+\Psi'(0)).
$$
So far they made no use of the important assumption $|\alpha|=|\beta|$. Coming to prove that $0<t<1$ this
assumption is needed. Indeed if $\alpha=r\cdot e^{i\theta}$ and $\beta=r\cdot e^{i\phi}$, where $0<\theta-\phi<2\pi$
(by $\alpha\ne\beta$ and $|\alpha|=|\beta|$) a simple computation gives
$$
\Psi'(0)=\pm\cos\frac{1}{2}(\theta-\phi).
$$
Hence $0<t<1$ and the elegant proof is done. \\
\\
Immediate consequences (see \cite{d}, Corollary 1 and Corollary 2 on page 287) are that each extreme point of $S$
and each support point of $S$ have the monotonic modulus property. We show how to get more information on $f$, based
on the above nice proof. The two functions $w\pm\Psi(w)$ are analytic and injective in $D$. In fact this is true
in every domain that is complementary to two disjoint slits that start respectively at $\alpha$ and at $\beta$
and extend to infinity. We note that if $w\not\in\{\alpha,\beta\}$ then also $w\pm\Psi(w)\not\in\{\alpha,\beta\}$
(for $w\pm\Psi(w)=\alpha\Rightarrow (w-\alpha)^2=(w-\alpha)(w-\beta)\Rightarrow w=\alpha\,\,{\rm or}\,\,w=\beta$).
Hence the following $4$ compositions are analytic, single-valued and injective in $D$ and omit $\{\alpha,\beta\}$,
$$
\begin{array}{l} g_1(w)=w+\Psi(w)+\Psi(w+\Psi(w)), \\ g_2(w)=w+\Psi(w)-\Psi(w+\Psi(w)), \\
g_3(w)=w-\Psi(w)+\Psi(w-\Psi(w)), \\ g_4(w)=w-\Psi(w)-\Psi(w-\Psi(w)).\end{array}
$$ 
These $4$ functions have disjoint images (for $\xi+\Psi(\xi)$ and $\eta-\Psi(\eta)$ are disjoint, so
$\eta=w_2+\Psi(w_2)$, $\xi=w_1+\Psi(w_1)$ give us the conclusion that $g_1, g_3$ are disjoint. Similarly
$\eta=w_2-\Psi(w_2)$, $\xi=w_1+\Psi(w_1)$ show that $g_1, g_4$ are disjoint. Also $\xi+\Psi(\xi)$ is injective
hence $g_2, g_3$ are disjoint because $w_1+\Psi(w_1)\ne w_2-\Psi(w_2)$, again because the disjointness
of the functions of Brickman and Wilken.). Clearly we have 
$$
\sum_{j=1}^4g_j(w)=4w.
$$
We define for $1\le j\le 4$ and $w\not\in\{\alpha,\beta\}$,
$$
h_j(w)=\frac{g_jw)-g_j(0)}{g'_j(0)},
$$
then
$$
\sum_{j=1}^4g'_j(0)\cdot h_j(w)=\sum_{j=1}^4g_j(w)-\sum_{j=1}^4g_j(0)=4w,
$$
and
$$
\sum_{j=1}^4\frac{1}{4}g'_j(0)=1\,\,\,\,\,{\rm by}\,\,\sum_{j=1}^4g'_j(w)=4.
$$
We conclude that if for $1\le j\le 4$ we have $g'_j(0)>0$ then 
$$
w=\sum_{j=1}^4\left(\frac{1}{4}g'_j(0)\right)h_j(w)
$$
is a strict convex combination (no zero coefficients) of the $h_j$, $1\le j\le 4$. Thus if $f\in S$
omits the values $\alpha$, $\beta$ so that $g'_j(0)>0$ for $1\le j\le 4$, then $f$ has the following
representation
$$
f=\alpha_1\cdot f_1+\alpha_2\cdot f_2+\alpha_3\cdot f_3+\alpha_4\cdot f_4,
$$
where $0<\alpha_j<1$, $\sum_{j=1}^4\alpha_j=1$ and $f_j$ are distinct functions in $S$ that omit
non-empty open sets. Here, as in Brickman and Wilken's proof, $f_j=g_j\circ f$, $1\le j\le 4$. So
we need to prove that $g'_j(0)>0$ for $1\le j\le 4$. We, once more, will make a use in the assumption 
$|\alpha|=|\beta|$ (which was already used by Brickman and Wilken in the first step of our iteration).
Let us compute $g'_j(0)$.
$$
g'_j(w)=1\pm\Psi'(w)\widehat{\mp}(1\pm\Psi'(w))\Psi'(w\pm\Psi(w)),
$$
where $\widehat{\mp}$ are signs not synchronized with the other sign changes in the expression.
$$
g'_j(w)=(1\pm\Psi'(w))\cdot(1\widehat{\mp}\Psi'(w\pm\Psi(w))),
$$
$$
g'_j(0)=(1\pm\Psi'(0))\cdot(1\widehat{\mp}\Psi'(\pm\Psi(0))).
$$
Now we have
$$
\Psi(0)=(\alpha\beta)^{1/2},\,\,\,\Psi'(0)=-\left(\frac{\alpha+\beta}{2(\alpha\beta)^{1/2}}\right),
\Psi'(\pm\Psi(0))=-\left(\frac{\alpha^{1/2}\mp\beta^{1/2}}{2\{\mp(\alpha\beta)^{1/2}\}^{1/2}}\right).
$$
Hence
$$
1\pm\Psi(0)=1\mp\left(\frac{\alpha+\beta}{2(\alpha\beta)^{1/2}}\right)=2(\Psi'(\pm\Psi(0)))^2,
$$
$$
g'_j(0)=2\Psi'(\pm\Psi(0))^2(1\widehat{\mp}\Psi'(\pm\Psi(0))).
$$
We denote $A=\Psi'(\pm\Psi(0))$ and then we need $0<2\cdot A^2\cdot(1\mp A)$. This happens when $-1<A<1$
and so $-1<\Psi'(\pm\Psi(0))<1$. This means that
$$
-1<\left(\frac{\alpha^{1/2}\mp\beta^{1/2}}{2\{\mp(\alpha\beta)^{1/2}\}^{1/2}}\right)<1,
$$
for if $\alpha=re^{i\theta}$, $b=re^{i\phi}$, then:
$$
\frac{e^{i\theta/2}\mp e^{i\phi/2}}{2\{\mp e^{i(\theta+\phi}/2\}^{1/2}}=\sin\left(\frac{\theta-\phi}{4}\right)\,\,
{\rm or}\,\,\cos\left(\frac{\theta-\phi}{4}\right).
$$
and we already know that this is indeed the case when $|\alpha|=|\beta|$. This proves the case $n=2$ in our general
theorem below.

\begin{theorem}\label{thm1}
If the function $f\in S$ omits two values of equal modulus, and if $n$ is a natural number, $n\in\mathbb{Z}^+$, then 
$f$ has the form
$$
f=\sum_{j=1}^{2^n}\alpha_j\cdot f_j,
$$
where $0<\alpha_j<1$, $1\le j\le 2^n$, $\sum_{j=1}^{2^n}\alpha_j=1$ and where $f_j$, $1\le j\le 2^n$ are different 
functions in $S$ that omit (each) open non-empty sets
\end{theorem}
\noindent
{\bf Proof.} \\
We denote $D=f(U)$ and we assume that $\alpha, \beta\not\in D$, $\alpha\ne\beta$, $|\alpha|=|\beta|$. We define
an analytic and single-valued function in $D$ by $\Psi(w)=\{(w-\alpha)(w-\beta)\}^{1/2}$ and denote two more functions
in $D$, $\Psi_1(w)=w+\Psi(w)$, $\Psi_2(w)=w-\Psi(w)$. We will define a sequence of $n$ sequences of functions. In
the $j$'th sequence there will be $2^j$ functions. The first sequence is: $g_{11}=\Psi_1$, $g_{12}=\Psi_2$. We now assume that
$j>1$ and that the $(j-1)$'st sequence is: $g_{(j-1)k}$, $k=1,2,\dots,2^{j-1}$. Then the $j$'th sequence is:
$$
\left\{\begin{array}{lll} g_{jk}(w) & = & (\Psi_1\circ g_{(j-1)k})(w) \\ g_{j(k+2^{j-1})}(w) & = & (\Psi_2\circ g(_{(j-1)k})(w)
\end{array}\right.\,\,\,1\le k\le 2^{j-1}.
$$
The functions in our series are injective ($g_{jk}(w_1)=g_{jk}(w_2)\Rightarrow (\Psi_l\circ g_{(j-1)k})(w_1)=
(\Psi_l\circ g_{(j-1)k})(w_2)$ for $l=1$ or $l=2$. $\Rightarrow g_{(j-1)k}(w_1)=g_{(j-1)k}(w_2)\Rightarrow
w_1=w_2$ inductively). The functions are pairwise disjoint in each of the $n$ sequences (functions within the same
sequence). For if $g_{jk}(w_1)=g_{jk}(w_2)$ then there can be only two possibilities: \\
(i) $(\Psi_s\circ g_{(j-1)k})(w_1)=(\Psi_s\circ g_{(j-1)l})(w_2)$ where $s=1$ or $s=2$. But $\Psi_l$ is injective
and hence $g_{(j-1)k}(w_1)=g_{(j-1)l}(w_2)$ and we use induction. \\
(ii) $(\Psi_1\circ g_{(j-1)k})(w_1)=(\Psi_2\circ g_{(j-1)l})(w_2)$ but $\Psi_1$ and $\Psi_2$ are disjoint and hence 
again $g_{(j-1)k}(w_1)=g_{(j-1)l}(w_2)$. \\
In particular for the $n$'th sequence we have: the functions $g_{nk}$, $1\le k\le 2^n$ are analytic, single-valued,
injective and disjoint in $D\subseteq\mathbb{C}-\{\alpha,\beta\}$. Also we have
$$
\sum_{k=1}^{2^n}g_{nk}(w)=2^n\cdot w.
$$
For we can use inductive argument as follows
$$
\sum_{k=1}^{2^n}g_{nk}(w)=\sum_{k=1}^{2^{n-1}}\{(\Psi_1\circ g_{(n-1)k})(w)+(\Psi_2\circ g_{(n-1)k})(w)\}=
$$
$$
=2\sum_{k=1}^{2^{n-1}}g_{(n-1)k}(w)=2\cdot(2^{n-1}\cdot w)=2^n\cdot w.
$$
We define for $1\le j\le 2^n$,
$$
h_j(w)=\frac{g_{nj}(w)-g_{nj}(0)}{g'_{nj}(0)},
$$
and then
$$
\sum_{j=1}^{2^n}g'_{nj}(0)\cdot h_j(w)=2^nw,\,\,{\rm and}\,\,\sum_{j=1}^{2^n}2^{-n}g'_{nj}(0)=1,
$$
where the second identity originates in
$$
\sum_{j=1}^{2^n}g'_{nj}(w)=(2^n\cdot w)'=2^n.
$$
We conclude that if for $1\le j\le 2^n$ we have $g'_{nj}(0)>0$, then
$$
w=\sum_{j=1}^{2^n}(2^{-n}g'_{nj}(0))\cdot h_j(w),
$$
the usual convex combination with positive coefficients of the $h_j(w)$'s. If this is the case, we define
$\alpha_j=2^{-n}g'_{nj}(0)$, $f_j=h_j\circ f$, $1\le j\le 2^n$ and we get the convex representation
$f=\sum_{j=1}^{2^n}\alpha_j\cdot f_j$ that we were looking for. For it is obvious that each $f_j\in S$
and those functions omit open non-empty sets, for the $h_j$ do, because the $g_{nj}$'s are disjoint. Thus we 
need to prove that $g'_{nj}(0)>0$ for $1\le j\le 2^n$. This follows by induction and by the assumption
that $|\alpha|=|\beta|$, $\alpha\ne\beta$. We note the following
$$
g'_{nj}(0)=\Psi_s'(g_{(n-1)j}(0))\cdot g'_{(n-1)j}(0)=\{1\pm\Psi'(g_{(n-1)j}(0))\}\cdot g'_{(n-1)j}(0)=
$$
$$
=\left\{1\pm\frac{2g_{(n-1)j}(0)-\alpha-\beta}{2\{(g_{(n-1)j}(0)-\alpha)(g_{(n-1)j}(0)-\beta)\}^{1/2}}\right\}
\cdot g'_{(n-1)j}(0)>0.
$$
We elaborate a bit more this final part of the proof. The proof that $g'_{nj}(0)>0$
is inductive (on $n$). It is convenient to denote $X_n=g_{nj}(0)$ ($j$ is fixed) and the induction assumption
is that $|X_n-\alpha|=|X_n-\beta|$. By $X_n=X_{n-1}\pm(X_n-\alpha)^{1/2}(X_n-\beta)^{1/2}$ we get
$$
\left\{\begin{array}{lll}
|X_n-\alpha| & = & |X_{n-1}-\alpha|^{1/2}|(X_{n-1}-\alpha)^{1/2}\pm(X_{n-1}-\beta)^{1/2}|, \\
|X_n-\beta| & = & |X_{n-1}-\beta|^{1/2}|(X_{n-1}-\alpha)^{1/2}\pm(X_{n-1}-\beta)^{1/2}|,
\end{array}\right.
$$
and hence $|X_n-\alpha|=|X_n-\beta|$ for all n. Hence
$$
\Psi'(g_{nj}(0))=\frac{2X_n-\alpha-\beta}{2\{(X_n-\alpha)(X_n-\beta)\}^{1/2}}=\frac{(X_n-\alpha)+(X_n-\beta)}
{2\{(X_n-\alpha)(X_n-\beta)\}^{1/2}},
$$
and we conclude that indeed $-1<\Psi'(g_{nj}(0))<1$. $\qed $ \\
\\
The construction in the proof of Theorem \ref{thm1} applies to any natural number $n\in\mathbb{Z}^+$.
A natural question is whether when $n\rightarrow\infty$ it converges to some kind of, say, an integral
representation of the function $f\in S$ that omits $\{\alpha,\beta\}$, where as usual $\alpha\ne\beta$,
$|\alpha|=|\beta|$. To start with, when we inquire if a recursion such as the one we have $g_{k+1}(w)=
w+\Psi(g_k(w))$ or $g_{k+1}(w)=w-\Psi(g_k(w))$ converges (the sign is chosen at each stage arbitrarily).
We first try to solve for $g$ in $g(w)=w+\Psi(g(w))$ or $g(w)=w-\Psi(g(w))$. We immediately note the
following,

\begin{proposition}\label{prop1}
Let us consider the following functions that result by applying finitely many times recursions of the
form
$$
g_0(w)=w,\,\,g_{k+1}(w)=w\pm\Psi(g_k(w)),\,\,k=0,1,2,\ldots.
$$
where at each step th sign $+$ or $-$ is chosen arbitrarily. Then all the resulting functions have
a unique fixed-point which is the same for all of them. This fixed-point is the rational function
$$
g(w)=\frac{w^2-\alpha\cdot\beta}{2w-\alpha-\beta}.
$$
\end{proposition}
\noindent
{\bf Proof.} \\
Solving for $g=w\pm\{(g-\alpha)(g-\beta)\}^{1/2}$ amounts in the equation $(g-w)^2=(g-\alpha)(g-\beta)$
regardless of the sign. This last equation is linear in $g$, $-2wg+w^2=-(\alpha+\beta)g+\alpha\beta$ and
it's (unique) solution is
$$
\frac{w^2-\alpha\cdot\beta}{2w-\alpha-\beta}.
$$The same is true when we solve the fixed-point equation of higher members of the recursion. For example,
solving for $g=w\pm\Psi(w\pm\Psi(g))$, is independent of the sign choices. It leads to
$$
(\alpha-\beta)^2(w^2-\alpha\beta)=(\alpha-\beta)^2(2w-\alpha-\beta)g.
$$
$\qed$ \\
\\
We conclude this section by noting that in passing with the sum of $2^n$ elements $\sum_{j=1}^{2^n}g'_{nj}(0)h_j(w)=
2^n\cdot w$ to the next sum, that of $2^{n+1}$ elements, $\sum_{j=1}^{2^{n+1}}g'_{(n+1)j}(0)\tilde{h}_j(w)=2^{n+1}\cdot w$, each
element in the former sum gave birth to two descendents $w+\Psi(g_{nj}(w))$ and $w-\Psi(g_{nj}(w))$. So in a sense,
each of the elements in a particular sum (say the one with $2^n$ elements) developed from a well-defined chain
of elements in the former (smaller) sums, in a way that resembles partial sums in s series development. When 
$n\rightarrow\infty$ we can interpret our recursive process as integrating all these multitude of elements that can be
thought of as the values of a  random variable over a probability space with the uniform distribution.

\section{One more property of support points of $S$}\label{sec3}

We recall that the space $H(U)$ is a linear topological locally convex space. The normalized conformal
mappings $S\subset H(U)$ is a compact topological subspace of $H(U)$. The topology is that of uniform convergence
on compact subsets. If $f\in S$ is a support point of $S$ that corresponds to the continuous linear functional
$L$ on $H(U)$, then by the definition $\Re L(f)=\max_{g\in S} \Re L(g)$. The complement of the image of
$f$, $\Gamma=\mathbb{C}-f(U)$ is an analytic curve having the property of increasing modulus and having the
$\pi/4$-property, i.e. the angle between the segment that connects the origin to the tip of $\Gamma$ is
at most $\pi/4$. Moreover, $\Gamma$ has an asymptotic direction at $\infty$. There is a point $z\in\partial U$ such
that $f$ is analytic on $\overline{U}-\{z_1\}$ and has a pole of order $2$ at $z=z_1$. Also, if $w_0$is the tip of
the slit $\Gamma$, then there is a point $z_0\in\partial U$such that $w_0=f(z_0)$, and $f'(z_0)=0$. If the functional $L$
is not constant on $S$ (as we assume throughout) then $L(f^2)\ne 0$ as is well known. In fact this was used 
in order to prove that the slit has an asymptotic direction at infinity. See \cite{s}. Also, in \cite{d}, Theorem 10.4
on page 307, and Theorem 10.5 on page 311.
It is here that we go further and prove a family of inequalities that involve $\Re L(\overline{z}_0^j (f(z)^{j+1})$.

\begin{theorem}\label{thm2}
Let $L$ be a continuous linear functional on $H(U)$ which is not constant on $S$. Let $f\in S$ satisfy
the equation $\Re L(f)=\max_{g\in S}\Re L(g)$, and suppose that $|z_0|=1$, $f'(z_0)=0$. Then for any natural number
$n\in\mathbb{Z}^+$, and for any positive real number $t\in\mathbb{R}^+$:
$$
\Re L\left\{\overline{z}_0f(z)^2e^{-t}\right\}+
\Re L\left\{\int_t^{\infty}\left\{\frac{e^s f(z,s)(k(s)f(z,s))^2}{1-k(s)f(z,s)}\right\}ds\right\}+o(e^{-t})\le 0.
$$
\end{theorem}
\noindent
{\bf Proof.} \\
Since $\Gamma=\mathbb{C}-f(U)$ is a slit, we can embed $f$ inside a L\"owner chain. We briefly recall
this standard procedure (see Chapter 3 in \cite{d}, 76-92). One chooses a parametric representation
of $\Gamma$, $w=\Psi(t)$, $0\le t<\infty$ so that $\Psi(0)=f(z_0)$, $\Psi(s)\ne\Psi(t)$ for $s\ne t$. Also,
if $\Gamma_t$ is the tail of $\Gamma$ from $\Psi(t)$ to $\infty$, then $g(z,t)$ is the Riemann mapping
of $U$ onto $\mathbb{C}-\Gamma_t$ so that $g(0,t)=0$, $g'(0,t)>0$ and we have:
$$
g(z,t)=e^t\left\{z+\sum_{n=2}^{\infty}b_n(t)z^n\right\},\,\,\,0\le t<\infty.
$$
We define 
$$
f(z,t)=g^{-1}(f(z),t)=e^{-t}\left\{z+\sum_{n=2}^{\infty}a_n(t)z^n\right\}.
$$
Then $f(z,t)$ is called a L\"owner chain and it satisfies:
$$
\frac{\partial f(z,t)}{\partial t}=-f(z,t)\cdot\frac{1+k(t)f(z,t)}{1-k(t)f(z,t)},
$$
$$
f(z,0)\equiv z,\,\,\,\forall\,z\in U,
$$
$$
\lim_{t\rightarrow\infty}e^tf(z,t)\equiv f(z),\,\,\,\forall\,z\in U,
$$
where the limit is uniform on compact subsets of $U$. The point $1/k(t)=\overline{k(t)}$ is that point on
$\partial U$ that is mapped by $f(z,t)$ onto the tip of $\Gamma_t$. We note that $e^tf(z,t)\in S$, $0\le t<\infty$
and so: \\
\\
(1) $\Re L(e^tf(z,t)-f(z))\le 0$, $0\le t<\infty$. \\
\\
(2) On the other hand we have:
$$
f(z)-e^tf(z,t)=\lim_{T\rightarrow\infty}\left\{e^Tf(z,T)-e^tf(z,t)\right\}=\left[e^sf(z,s)\right]_{s=t}^{\infty}=
\int_t^{\infty}h(z,s)ds,
$$
where $\int h(z,s)ds=e^sf(z,s)$. \\
\\
(3) By differentiation:
$$
h(z,s)=\frac{\partial}{\partial s}\left\{e^s f(z,s)\right\}=e^s f(z,s)+e^s\frac{\partial f(z,s)}{\partial s}=
$$
$$
=e^s f(z,s)-e^s f(z,s)\frac{1+k(s)f(z,s)}{1-k(s)f(z,s)}=-e^s f(z,s)\frac{2k(s)f(z,s)}{1-k(s)f(z,s)}.
$$
(4) From the equations in (1), (2) and (3) we conclude that:
$$
\Re L\left(\int_t^{\infty}\left\{e^s f(z,s)\frac{k(s)f(z,s)}{1-k(s)f(z,s)}\right\}ds\right)\le 0.
$$
We recall that $\lim_{s\rightarrow\infty}k(s)=\overline{z}_0$ and also $\lim_{s\rightarrow\infty}e^s f(z,s)=f(z)$
uniformly on compact subsets of $U$. So we can write $e^sf(z,s)=f(z)+\epsilon(s)$, where $\lim_{s\rightarrow\infty}\epsilon(s)=0$.
Also $k(s)=\overline{z}_0+\delta(s)$, where $\lim_{s\rightarrow\infty}\delta(s)=0$. Hence:
$$
e^s f(z,s)\frac{k(s)f(z,s)}{1-k(s)(\overline{z}_0+\delta(s))e^{-s}(f(z)+\epsilon(s))f(z,s)}=
$$
$$
=(f(z)+\epsilon(s))
\frac{(\overline{z}_0+\delta(s))e^{-s}(f(z)+\epsilon(s))}{1-(\overline{z}_0+\delta(s))e^{-s}(f(z)+\epsilon(s))}=
$$
$$
=\overline{z}_0f(z)^2e^{-s}+\left(\overline{z}_0(2\epsilon(s)f(z)+\epsilon(s)^2)+\delta(f(z)+\epsilon(s))^2\right)e^{-s}+
$$
$$
+\frac{e^s f(z,s)(k(s)f(z,s))^2}{1-k(s)f(z,s)}.
$$
Integrating between $t$ and $\infty$ we obtain the identity:
$$
\int_t^{\infty}\left\{e^s f(z,s)\frac{k(s)f(z,s)}{1-k(s)f(z,s)}\right\}ds=\overline{z}_0f(z)^2e^{-t}+
\int_t^{\infty}\left\{\frac{e^s f(z,s)(k(s)f(z,s))^2}{1-k(s)f(z,s)}\right\}ds+
$$
$$
+\int_t^{\infty}\left(\overline{z}_0(2\epsilon(s)f(z)+\epsilon(s)^2)+\delta(f(z)+\epsilon(s))^2\right)e^{-s}ds.
$$
The last integral is $o(e^{-t})$ for $t\rightarrow\infty$. Hence using the inequality in (4), we obtain:
$$
\Re L\left\{\overline{z}_0f(z)^2e^{-t}\right\}+
\Re L\left\{\int_t^{\infty}\left\{\frac{e^s f(z,s)(k(s)f(z,s))^2}{1-k(s)f(z,s)}\right\}ds\right\}+o(e^{-t})\le 0.
$$
This proves the theorem. $\qed $ \\

\begin{remark}\label{rem1}
It is not possible to deduce from the inequality in Theorem \ref{thm2} that 
$$
\Re L\left\{\overline{z}_0f(z)^2e^{-t}\right\}+o(e^{-t})\le 0,
$$
when $t\rightarrow\infty$ as the author wrongly thought in the first version of this
paper. The author thanks the referee for his remark and insight on that matter. Here is a simple
example that shows that such an inequality can not be true. Let $L(a_0+a_1z+a_2z^2+\ldots)$ be
the third coefficient functional on $H(U)$. Then as is well known, $L$ has two support points
in $S$, $f_1(z)=z/(1-z)^2$ the Koebe function and its rotation $f_2(z)=-f_1(-z)=z/(1+z)^2$. For $f_1(z)$
we have $z_0=-1$ ($f_1'(-1)=0$), and for $f_2(z)$, $z_0=1$. Clearly $L(f_1^2)=L(f_2^2)=1$. Hence
$\Re L(\overline{z}_0f_1^2)=-1\cdot 1=-1<0$, but $\Re L(\overline{z}_0f_2^2)=1\cdot 1=1>0$.
\end{remark}

\section{Properties of support points are inherited by less geometric families of mappings}\label{sec4}

If $f\in S$is a support point of $f$ that corresponds to the continuous linear functional $L$ on $H(U)$,
then $\Gamma=\mathbb{C}-f(U)$is an analytic curve, called the slit of $f$. It starts at its tip $w_0$ and
monotonically extends to infinity. The tip $w_0$ has a single pre-image $z_0$ on $\partial U$, $w_0=f(z_0)$
and $f'(z_0)=0$. In this section we will see that the tip of the slit, $w_0$, already appears in a support point
of univalent polynomials.

\begin{definition}\label{def2}
Let $n\in\mathbb{Z}^+$ be a natural number. The family of all the polynomials in $S$, of degree $n$ or less
will be denoted by $S_n$. We note that $S_n$ is a compact subspace of $H(U)$.
\end{definition}

\begin{example}\label{ex1}
To demonstrate Definition \ref{def2} we note that $S_1=\{z\}$ and $S_2=\{z+\alpha z^2\,|\,|\alpha|\le 1/2\}$.
\end{example}
\noindent
Here is our result.

\begin{theorem}\label{thm4}
Let $n>1$ be a natural number. Let $L\in H(U)'$ be a continuous linear functional. Then either $L$ is constant
on $S_n$, or, if $L$ is not constant on $S_n$, then if $p\in S_n$ solves the following
extremal problem $|L(p)|=\max_{f\in S_n}|L(f)|$, then $p'(z)$ has a zero on $\partial U$.
\end{theorem}
\noindent
{\bf Proof.} \\
Let us suppose that $L$ is not constant on $S_n$. Then there are two polynomials $q_2, q_1\in S_n$ such that
$q_2-q_1\not\in\mathbb{C}$ and $L(q_2)\ne L(q_1)$. Using the normalization of elements in $S_n$ we have
$q_1(0)=q_2(0)=0$ and $q_1'(0)=q_2'(0)=1$. Hence $q(0)=q'(0)=0$ and $L(q)=L(q_2-q_1)=L(q_2)-L(q_1)\ne 0$.
We proceed with a Rouche's type of principle for injectivity.
\begin{lemma}\label{lem1}
If $f(z)=z+a_2 z^2+\ldots\in S$ is analytic in a neighborhood of $\overline{U}$ such that $f'(z)$ does not
have zero on $\partial U$, then for any function $g(z)$ analytic in a neighborhood of $\overline{U}$ there
exists a $\delta>0$ (depending on $g$), so that if $|w_0|<\delta$, then $f(z)+w_0\cdot g(z)$ is injective
on $\overline{U}$.
\end{lemma}
\noindent
{\bf A proof of Lemma \ref{lem1}.} \\
Let us fix $w_0$. We denote $F(z)=f(z)+w_0\cdot g(z)$. Then for any $z, w\in\overline{U}$ we have:
$$
|F(z)-F(w)|=|f(z)-f(w)+w_0(g(z)-g(w))|=
$$
$$
=|f(z)-f(w)|\times\left|1+w_0\left(\frac{g(z)-g(w)}{f(z)-f(w)}\right)\right|,
$$
where for $z=w$ we agree to interpret
$$
\frac{g(z)-g(w)}{f(z)-f(w)}=\frac{g'(z)}{f'(z)}.
$$
Since $f$ is in fact injective on $\overline{U}$ (because $f'(z)$ does not vanish on $U$) we deduce that
for $z\ne w$, $f(z)-f(w)\ne 0$ and so it is sufficient to prove that
$$
\frac{g(z)-g(w)}{f(z)-f(w)},
$$
is bounded on $\overline{U}\times\overline{U}$. We write the following identity:
$$
\frac{g(z)-g(w)}{f(z)-f(w)}=\left(\frac{g(z)-g(w)}{z-w}\right)\left/\left(\frac{f(z)-f(w)}{z-w}\right)\right..
$$
So it is sufficient to prove that:
$$
\max_{\overline{U}\times\overline{U}}\left|\frac{g(z)-g(w)}{z-w}\right|\le M<\infty
$$
and
$$
\min_{\overline{U}\times\overline{U}}\left|\frac{f(z)-f(w)}{z-w}\right|\ge\epsilon>0.
$$
For the minimum. If the estimate is false then there is a sequence $(z_k,w_k)\in\overline{U}\times\overline{U}$
so that
$$
\frac{f(z_k)-f(w_k)}{z_k-w_k}\rightarrow 0.
$$
The set $\overline{U}\times\overline{U}$ is compact in $\mathbb{C}\times\mathbb{C}$ and so we may assume
that $z_k\rightarrow a$ and $w_k\rightarrow b$ and we get in the case $a\ne b$ the equation
$$
\frac{f(a)-f(b)}{a-b}=0,
$$
which contradicts the injectivity of $f$ in $\overline{U}$. If $a=b$ we get a contradiction to the assumption
that $f'(z)$ does not have a zero in $\overline{U}$. For the maximum. Using arguments similar to those
above we get $a, b\in\overline{U}$ such that 
$$
\left|\frac{g(a)-g(b)}{a-b}\right|=\infty.
$$
If $a\ne b$ this contradicts the fact that $g$ is analytic in a neighborhood of $\overline{U}$. If $a=b$ this
contradicts the fact that $g'$ is analytic in a neighborhood of $\overline{U}$. The proof of Lemma \ref{lem1}
is now completed. $\qed $ \\
\\
We now conclude the proof of Theorem \ref{thm4} as follows. Using Lemma \ref{lem1} there exists an $\epsilon>0$
such that $p(z)+\epsilon e^{i\theta} q(z)\in S_n$, for any $0\le\theta<2\pi$. We use here the assumption that
$p'(z)$ does not have a zero on $\overline{U}$. Since $p$ is solving the extremal problem for $L$ we conclude
that $|L(p)|+\epsilon|L(q)|\le|L(p)|$ (one needs to choose properly the $\theta$). By $\epsilon>0$ it follows
that $L(q)=0$. This contradicts the assumption on $q$. Hence $p'(z)$ must have a zero on $\partial U$. $\qed $ \\

\begin{remark}\label{rem2}
(1) In particular for the coefficients functionals, $L(p)=p^{(j)}(0)/j!$ where $j$ is in $1<j\le n$ we note
that $L(z^j)=1\ne 0$, so we can use in Theorem \ref{thm4} $q(z)=z^j$. We conclude that if $p(z)$ maximizes
$|p^{(j)}(0)/j!|$ then $p'(z)$ must have a zero on $\partial U$. \\
(2) If, as in Theorem \ref{thm4}, $p'(e^{i\theta})=0$, then necessarily $p''(e^{i\theta})\ne 0$ for $p\in S_n$.
Thus $p'$ has only simple zeros on $\partial U$ and at least one. The reason that $p''(e^{i\theta})\ne 0$ for $p\in S_n$,
is that otherwise $p(z)$ is locally at $e^{i\theta}$ equivalent to $(z-e^{i\theta})^3$. This means that it
folds the tangent line to $\partial U$ at $e^{i\theta}$ by more than $\pi$ radians hence it can not be univalent
in a $U$-neighborhood of $e^{i\theta}$.
\end{remark}

\noindent
{\it Ronen Peretz \\
Department of Mathematics \\ Ben Gurion University of the Negev \\
Beer-Sheva , 84105 \\ Israel \\ E-mail: ronenp@math.bgu.ac.il} \\ 
 
\end{document}